\documentclass[12pt]{article}
\usepackage[english]{babel}
\usepackage{amsmath,amsthm}
\usepackage{amsfonts}
\usepackage{amssymb}

\theoremstyle{definition}

\theoremstyle{remark}

\numberwithin{equation}{section}

\begin{document}

\begin{center}
\bfseries  About the Kronecker’s Theorem  \\
\end{center}

\begin{center}

\vspace{5mm}
JinHua  Fei\\
\vspace{5mm}

ChangLing  Company of Electronic Technology    \, Baoji \,  Shannxi  \,  P.R.China \\
\vspace{5mm}

E-mail: feijinhuayoujian@msn.com \\

\end{center}

\vspace{3mm}
{\bfseries Abstract.}\,In this paper, we give a counter-example, in the general case, Kronecker theorem will derive contradiction. Kronecker’s theorem be correct after removing some conditions.


{\bfseries MR(2010)  Subject  Classification \quad 11J13 \, 11M06 } \\

\vspace{8mm}

In the references [1], Tom M. Apostol introduced the following famous theorem.

{\bfseries Kronecker’s theorem.}  If $ \alpha_1,\ldots,\alpha_n  $   are arbitrary real numbers, if  $ \theta_1,\ldots, \theta_n $ are linearly independent real numbers, and if $ \epsilon>0  $ is arbitrary, then there exists a real number $t$ and integers $ h_1,\ldots, h_n   $ such that

$$  \left|  t\theta_i -h_i - \alpha_i   \right | < \varepsilon    \qquad  for \,\,  i=1,2,\ldots,n $$

See page 150 of references [1].\\

We give a counter-example, in the general case, Kronecker theorem will derive contradiction.The theorem be true if the real numbers $ \alpha_1,\ldots, \alpha_n   $ are all zeros, for this question, the reader can refer to the page 541 of the references [3]. \\

For the sake of narration, let's give a few simple lemmas.\\

{\bfseries Lemma 1.} (1) When $ Res>1  $,  we have

$$ \log \zeta (s) = \sum_{n=2}^\infty  \frac{ \Lambda(n)  }{ n^s \log n  }  \qquad   \qquad   \qquad    \qquad  $$

Where $ \zeta(s) $ is Riemann zeta function and $ \Lambda(n) $  is von Mangoldt function.

See page23 of references [2].\\

(2) The inequalities

$$ \frac{1}{\sigma-1} < \zeta (\sigma) <  \frac{ \sigma}{ \sigma-1 }  \qquad   \qquad   \qquad    \qquad   $$

hold for all $ \sigma>0 $  . See page25 of references [2]. \\

{\bfseries Lemma 2.}  (1) When $ x \geq 2 $ , there is a constant $ c_1 > 0 $ such that

$$  \sum_{ 2\leq p \leq x }  \log p = x + O\left( x \exp (-c_1 \sqrt{\log x })   \right)   \qquad    $$

Where $ p $ is the prime number. See page179 of references [2].\\

(2) When $ x \geq 7 $  , there is a constant $ c_2 > 0 $  such that

$$  \sum_{ 2\leq p \leq x }  \frac{1}{p} = \log\log  x + b + O\left( \exp (-c_2 \sqrt{\log x })   \right) $$

where $ b $ is a constants. See page182 of references [2]. \\

{\bfseries Lemma 3.} (1) When $ 0< \delta \leq \frac{1}{20}  $  and    $  \log y =\frac{1}{\delta}  $  , we have

$$ \sum_ { y < p } \frac{1}{p^{1+\delta}}  \ll  1     \qquad   \qquad   \qquad    \qquad   \qquad  \qquad   $$

(2)  we have

$$  \sum_{ k= 2}^\infty \frac{1}{k} \sum_{2 \leq p} \frac{1}{p^k} \ll 1    \qquad   \qquad   \qquad    \qquad  \qquad                   $$  \\

Proof of (1).  By (1) of lemma 2, we have

$$ \sum_ { y < p } \frac{1}{p^{1+\delta}} = \int_y^\infty  \frac{ u^{-1-\delta}  }{ \log u  }  d \left( \sum_{p\leq u} \log p \right)  = - \frac{ y^{-1-\delta}}{ \log y  } \left( \sum_{p\leq y} \log p \right) $$

$$ - \int_y^\infty \left(  (-1-\delta) \frac{ u^{-2-\delta} }{ \log u  } - \frac{ u^{-2-\delta} }{ \log^2 u  } \right) \left( \sum_{p\leq u} \log p \right) du  \quad    \qquad   $$

 $$  \ll  \frac{ y^{-\delta} }{ \log y } + \int_y^\infty  \frac{ u^{-1-\delta}  }{ \log u } du \ll  \frac{ y^{-\delta} }{ \delta \log y } \ll 1  \qquad   \qquad   \qquad \quad $$

this completes the proof .\\

Proof of (2).  By (1) of lemma 2 and when $ k\geq 2 $, we have

$$  \sum_{2\leq p} \frac{1}{ p^k } \ll  \sum_{2\leq p} \frac{\log p}{ p^k } = \int_2^\infty u^{-k} d \left( \sum_{p\leq u}  \log p \right)  \quad $$

$$ \leq k \int_2^\infty u^{-1-k}   \left( \sum_{p\leq u}  \log p \right) du   \ll k \int_2^\infty  u^{-k } du     $$

$$   = \frac{k}{1-k}  \int_2^\infty d u^{1-k } = \frac{k}{k-1} 2^{1-k} \ll 2^{-k}   \qquad  $$

therefore

$$  \sum_{ k= 2}^\infty \frac{1}{k} \sum_{2 \leq p} \frac{1}{p^k} \ll   \sum_{ k= 2}^\infty \frac{1}{k} \,  2^{-k} \ll 1   \qquad  \qquad \qquad      $$

this completes the proof. \\

{\bfseries Proof of the Conclusion.}

Let $ \tau $ be  any real number, $ 0<\delta \leq \frac{1}{20} $ , by (1) of lemma 1, we have

$$  \zeta (1+\delta -i \tau)  = \sum_{n=2}^\infty \frac{ \Lambda (n) }{ n^{1+\delta -i \tau} \log n }   \qquad  \qquad \qquad \qquad \quad   $$

$$ =   \sum_{n=2}^\infty \frac{ \Lambda (n) }{ n^{1+\delta} \log n } \cos ( \tau \log n )  + i \sum_{n=2}^\infty \frac{ \Lambda (n) }{ n^{1+\delta} \log n } \sin ( \tau \log n )  $$

we take the imaginary part of both sides of this equation, then

$$  \arg   \zeta (1+\delta -i \tau)   = \sum_{n=2}^\infty \frac{ \Lambda (n) }{ n^{1+\delta} \log n } \sin ( \tau \log n ) $$

By lemma 3 and $ \log y = \frac{1}{\delta}  $ , we have

$$ \sum_{n=2}^\infty \frac{ \Lambda (n) }{ n^{1+\delta} \log n } \sin ( \tau \log n ) = \sum_{2\leq p} \frac{1}{p^{1+\delta}}  \sin ( \tau \log p ) \qquad  \quad $$

$$  +  \sum_{k=2}^\infty  \frac{1}{k} \sum_{2\leq p} \frac{1}{p^{k(1+\delta)}}  \sin ( \tau k \log p )    \qquad  \qquad \qquad  \quad   \qquad  \quad    $$

$$ = \sum_{2\leq p \leq y} \frac{1}{p^{1+\delta}}  \sin ( \tau \log p )  +   \sum_{ y < p } \frac{1}{p^{1+\delta}}  \sin ( \tau \log p ) + O(1) $$

$$ = \sum_{2\leq p \leq y} \frac{1}{p^{1+\delta}}  \sin ( \tau \log p )  +  O(1)  \qquad  \qquad \qquad \qquad \quad    $$

Easily seen

$$  \sin(\tau \log p) =  \sin \left(\tau \log p - \frac{\pi}{2} +\frac{\pi}{2} \right )    \qquad  \qquad \qquad    $$

$$ =  \cos \left(\tau \log p - \frac{\pi}{2} \right ) = \cos \left( 2\pi \left( \frac{1}{2\pi} \tau \log p - \frac{1}{4}  \right) \right )  $$

therefore

$$  \sum_{2\leq p \leq y} \frac{1}{p^{1+\delta}}  \sin ( \tau \log p ) = \sum_{2\leq p \leq y} \frac{1}{p^{1+\delta}}  \cos \left( 2\pi \left( \frac{1}{2\pi} \tau \log p - \frac{1}{4}  \right) \right )     $$

Evident $ \frac{1}{2\pi}  \log p \,\, ( 2\leq p \leq y ) $    are linearly independent over the integers. \\

By Kronecker’s theorem, there is the real number $ \tau_0 $  and the integer  $ h(p) $, we have

$$ \left|  \frac{1}{2\pi} \tau_0 \log p - \frac{1}{4} - h(p) \right| < \delta \,\,\,\,\, (2\leq p  \leq y ) $$ \\

Now, we take $ \tau=\tau_0 $ . Because

$$  \cos \left( 2\pi \left( \frac{1}{2\pi} \tau_0  \log p - \frac{1}{4}  \right) \right ) =   \cos \left( 2\pi \left( \frac{1}{2\pi} \tau_0  \log p - \frac{1}{4}  - h(p) \right) \right )  $$

$$  =   \cos \left( 2\pi \left| \frac{1}{2\pi} \tau_0  \log p - \frac{1}{4}  - h(p) \right| \right )  = 1 + O(\delta^2)  \qquad  \qquad \qquad   $$

therefore

$$  \sum_{2\leq  p \leq  y} \frac{1}{p^{ 1+\delta }} \cos \left( 2\pi \left( \frac{1}{2\pi} \tau_0  \log p - \frac{1}{4}  \right) \right )   \qquad  \qquad \qquad   $$

$$ = \sum_{2\leq  p \leq  y} \frac{1}{p^{ 1+\delta }} +O \left( \delta^2 \sum_{2\leq  p \leq  y} \frac{1}{p^{ 1+\delta }} \right )   \qquad  \qquad \qquad \quad  $$ \\

By (2) of lemma 2 and $ \log y = \frac{1}{\delta}  $ , we have

$$  \delta^2 \sum_{2\leq  p \leq  y} \frac{1}{p^{ 1+\delta }} \leq \delta^2 \sum_{2\leq  p \leq  y} \frac{1}{p}  \ll \delta^2 \log\log y \ll  \delta^2 \log \delta^{-1 }  \ll 1  $$ \\

By lemma3 and (1) of lemma1, we have

$$  \sum_{2\leq  p \leq  y} \frac{1}{p^{ 1+\delta }} =  \sum_{2\leq  p } \frac{1}{p^{ 1+\delta }} - \sum_{y <  p } \frac{1}{p^{ 1+\delta }} = \sum_{2\leq  p } \frac{1}{p^{ 1+\delta }}  + O(1)  \qquad  \qquad  \quad    $$

$$  =   \sum_{ n=2 }^\infty \frac{ \Lambda (n) }{ n^{1+\delta} \log n } - \sum_{k=2}^\infty  \frac{1}{k} \sum_{2\leq p} \frac{1}{p^{k(1+\delta)} } + O(1)  = \log \zeta (1+\delta) + O(1) $$ \\

From the above argument, we have

$$  \arg \zeta (1+\delta -i \tau_0) =  \log \zeta (1+\delta) + O(1)    \qquad    $$

because

$$\arg \zeta (1+\delta -i \tau_0) \leq 2\pi    \qquad  \qquad   \qquad  \qquad     \quad    $$

therefore

$$  c_3 \geq  \log  \zeta (1+\delta)   \qquad  \qquad   \qquad  \qquad   \qquad  \quad   $$

By (2) of lemma1, we have

$$ e^{c_3} \geq    \zeta (1+\delta)  \geq  \frac{1}{\delta}  \qquad  \qquad   \qquad  \qquad   \qquad     $$

when  $ \delta $  is sufficiently small, there is a contradiction. This completes the proof of the conclusion.

\vspace{10mm} \centerline{ REFERENCES } \vspace{5mm}

[1]  Tom M. Apostol,{ \itshape  Modular Function and Dirichlet Series in Number Theory, } Springer-Verlag New York, 1990 \\

[2]  Hugh L. Montgomery, Robert C. Vaughan,{ \itshape Multiplicative Number Theory I. Classical Theory,} Cambridge University Press, 2006. \\

[3]  Hua Loo Keng,{  \itshape  Introduction to Number Theory,} Springer-Verlag Berlin Heidelberg New York, 1982. \\

\end{document}